\documentclass[a4paper,11pt]{scrartcl}
\usepackage{graphicx}
\usepackage{amsmath,amssymb}
\usepackage{amsthm}
\usepackage{enumitem}
\usepackage{algpseudocode,algorithm}
\usepackage{graphicx}
\graphicspath{{./fig/}}
\usepackage[square,sort&compress,numbers]{natbib}
\usepackage{mathtools}\mathtoolsset{showonlyrefs=false}
\usepackage{url}

\usepackage{subfigure}

\newcommand{\Authornote}{\renewcommand{\thefootnote}{\fnsymbol{footnote}}}
\newcommand{\authornote}{\Authornote\footnote}
\newcommand{\refalg}[1]{Algorithm~\ref{#1}}

\newcommand{\reffig}[1]{Figure~\ref{#1}}

\newcounter{alnum}

\newenvironment{Problem}{\begin{array}.{*{20}{l}}\}}{\end{array}}

\newcommand{\MIN}{\mathop{\mathrm{Minimize}}}
\newcommand{\MAX}{\mathop{\mathrm{Maximize}}}

\newcommand{\ST}{\mathop{\mathrm{subject~to}}}

\newcommand{\diag}{\mathop{\mathrm{diag}}\nolimits}

\newcommand{\argmin}{\operatornamewithlimits{\mathrm{arg\,min}}}

\renewcommand{\Re}{\ensuremath{\mathbb{R}}}

\newcommand{\bi}[1]{\ensuremath{\boldsymbol{#1}}}

\newcommand{\pdif}[2]{\frac{\partial #1}{\partial #2}}

\begin{document}

\begin{center}
  {\Large\bfseries\sffamily%
  An Accelerated Uzawa Method for Application 
  }
  \par\medskip
  {\Large\bfseries\sffamily%
  to Frictionless Contact Problem 
  }%
  \par\bigskip
  Yoshihiro Kanno \authornote[2]{
    Mathematics and Informatics Center, 
    The University of Tokyo, 
    Hongo 7-3-1, Tokyo 113-8656, Japan.
    E-mail: \texttt{kanno@mist.i.u-tokyo.ac.jp}. 
  }
\end{center}

\begin{abstract}
  The Uzawa method is a method for solving constrained optimization 
  problems, and is often used in computational contact mechanics. 
  The simplicity of this method is an advantage, but its convergence is 
  slow. 
  This paper presents an accelerated variant of the Uzawa method. 
  The proposed method can be viewed as application of an accelerated 
  projected gradient method to the Lagrangian dual problem. 
  Preliminary numerical experiments suggest that the convergence of the 
  proposed method is much faster than the original Uzawa method. 
\end{abstract}

\begin{quote}
  \textbf{Keywords}
  \par
  Fast first-order method, 
  accelerated gradient scheme, 
  contact mechanics, 
  nonsmooth mechanics, 
  Uzawa method, 
  convex optimization. 
\end{quote}

\section{Introduction}

It has been recognized well that contact mechanics has close relation 
with optimization and variational inequalities \citep{DL76,Wri06}. 
The static frictionless contact problem of a linear elastic body, also 
called Signorini's problem, is one of the most fundamental problems in 
contact mechanics. 
This is a boundary value problem to find the equilibrium configuration 
of an elastic body, where some portion of the boundary of the body can 
possibly touch the surface of a rigid obstacle (or the surface of 
another elastic body). 
Positive distance between the elastic body and the obstacle surface 
(i.e., positive {\em gap\/}) implies zero contact pressure (i.e., zero 
{\em reaction\/}), while nonzero reaction implies zero gap. 
This disjunction nature can be described by using complementarity conditions. 
Moreover, the frictionless contact problem can be formulated as a 
continuous optimization problem under inequality constraints \citep{Wri06}. 

The Uzawa method is known as a classical method for solving constrained 
optimization problems \citep{All07,Cia89,Uza58}. 
Due to ease in implementation, the Uzawa method is often applied to 
contact problems \citep{TWH12,PLS08,Kok09,HKD04,Tem12,VN07,Rud15,Rao99}. 
Major drawback of the Uzawa method is that its convergence is slow; 
it exhibits only linear convergence in general. 

Recently, accelerated, or ``optimal'' \citep{Nes83,Nes04}, first-order methods 
have received substantial attention, particularly for solving 
large-scale optimization problems; see, e.g., 
\cite{BT09,BBC11,OdC15,LS13}. 
Advantages of most of these methods include 
ease of implementation, 
cheap computation per each iteration, 
and fast local convergence. 
Application of an accelerated first-order method to computational 
mechanics can be found in \cite{Kan16}. 
In this paper, we apply the acceleration scheme in \cite{BT09} 
to the Uzawa method.


The paper is organized as follows. 
Section~\ref{sec:fundamental} provides an overview of the necessary
background of the frictionless contact problem and the Uzawa method. 
Section~\ref{sec:acceleration} presents an accelerated Uzawa method. 
Section~\ref{sec:ex} reports the results of preliminary numerical 
experiments. 
Some conclusions are drawn in section~\ref{sec:conclude}.


In our notation, 
${}^{\top}$ denotes the transpose of a vector or a matrix. 
For two vectors $\bi{x}=(x_{i}) \in \Re^{n}$ and 
$\bi{y} = (y_{i}) \in \Re^{n}$, we define 
\begin{align*}
  \min \{\bi{x}, \bi{y}\} 
  = (\min \{ x_{1},y_{1} \}, \dots, \min \{ x_{n},y_{n} \})^{\top} . 
\end{align*}
We use $\Re_{-}^{n}$ to denote the nonpositive orthant, i.e., 
$\Re_{-}^{n} = \{ \bi{x} = (x_{i}) \in \Re^{n} \mid x_{i} \le 0 \ (i=1,\dots,n) \}$.
For $\bi{y} \in \Re^{n}$, we use $\Pi_{-}(\bi{y})$ to denote the 
projection of $\bi{y}$ on $\Re_{-}^{n}$, i.e., 
\begin{align*}
  \Pi_{-}(\bi{y})
  = \argmin \{ \| \bi{x} - \bi{y} \| \mid \bi{x} \in \Re_{-}^{n} \} 
  = \min \{ \bi{y}, \bi{0} \} .
\end{align*}
We use $\diag(\bi{x})$ to denote a diagonal matrix, the vector of 
diagonal components of which is $\bi{x}$.

\section{Fundamentals of frictionless contact and Uzawa method}
\label{sec:fundamental}

We briefly introduce the frictionless contact problem of an elastic body; 
see, e.g., \cite{Wri06} and \cite{DL76} for fundamentals of contact mechanics. 

Consider an elastic body subjected to a static load and a rigid obstacle 
fixed in space. 
Since the body cannot penetrate the surface of the obstacle, the 
deformation of the body is constrained from one side by the obstacle 
surface. 
We also assume the absence of friction and adhesion between the body 
surface and the obstacle surface. 
The set of these conditions is called the frictionless unilateral 
contact. 

Suppose that the conventional finite element procedure is adopted for 
discretization of the elastic body. 
Let $\bi{u} \in \Re^{d}$ denote the nodal displacement vector, where $d$ 
is the number of degrees of freedom of displacements. 
We use $\pi(\bi{u})$ to denote the total potential energy caused by $\bi{u}$. 
At the (unknown) equilibrium state, the nodes on a portion of the body 
surface can possibly make contact with the obstacle surface. 
Such nodes are called the contact candidate nodes, and the number of 
them is denoted by $m$. 
We use $g_{i}(\bi{u})$ $(i=1,\dots,m)$ to denote the gap between 
the $i$th contact candidate node and the obstacle surface. 
The non-penetration conditions are then formulated as 
\begin{align}
  g_{i}(\bi{u}) \ge 0 , 
  \quad  i=1,\dots,m .
  \label{eq.gap.1}
\end{align}
The displacement vector at the equilibrium state minimizes the total 
potential energy under the constraints in \eqref{eq.gap.1}. 
Namely, the frictionless contact problem is formulated as follows: 
\begin{subequations}\label{P.contact.1}%
  \begin{alignat}{3}
    &  \MIN &{\quad}& 
    \pi(\bi{u}) \\
    & \ST && 
    g_{i}(\bi{u}) \ge 0 , 
    \quad i=1,\dots,m . 
  \end{alignat}
\end{subequations}

The Uzawa method solving problem \eqref{P.contact.1} is listed in 
\refalg{alg:contact.1} \citep{All07,Cia89,Uza58,AC91}. 
It is worth noting that step~\ref{alg:contact.1.equilibrium} of 
\refalg{alg:contact.1} can be performed as the equilibrium analysis with 
a conventional finite element code. 
Also, step~\ref{alg:contact.1.gradient} is simple to implement. 
Due to such ease in implementation, the Uzawa method is widely used 
in computational contact mechanics. 

\begin{algorithm}
  \caption{Uzawa method}
  \label{alg:contact.1}
  \begin{algorithmic}[1]
    \Require
    $\alpha > 0$, $\bi{r}^{(0)} \in \Re_{-}^{m}$
    \For{$k=0,1,\dots$}
    \State \label{alg:contact.1.equilibrium}
    $\bi{u}^{(k)}$ solves 
    $\displaystyle 
    \nabla \pi(\bi{u}) + \sum_{i=1}^{m} r^{(k)}_{i} \nabla g_{i}(\bi{u}) = \bi{0}$
    \State \label{alg:contact.1.gradient}
    $r_{i}^{(k+1)} 
    := \min \{ 0, r_{i}^{(k)} + \alpha g_{i}(\bi{u}^{(k)}) \}$ 
    $(i=1,\dots,m)$
    \EndFor
  \end{algorithmic}
\end{algorithm}

As shown in \citet[section~9.4]{Cia89}, the Uzawa method can be viewed as 
a projected gradient method solving the Lagrange dual problem. 
Essentials of this observation are repeated here. 
The Lagrangian $L : \Re^{d} \times \Re^{m} \to \Re \cup \{ -\infty \}$ 
associated with problem \eqref{P.contact.1} is given by 
\begin{align}
  L(\bi{u}, \bi{r}) = 
  \begin{dcases*}
    \pi(\bi{u}) + \sum_{i=1}^{m} r_{i} g_{i}(\bi{u}) 
    & if $\bi{r} \le \bi{0}$, \\
    -\infty 
     & otherwise. 
  \end{dcases*}
  \label{eq.def.Lagrangian.1}
\end{align}
Here, the Lagrange multipliers, $r_{1},\dots,r_{m}$, correspond to the 
reactions. 
The inequality constraints imposed on the reactions in 
\eqref{eq.def.Lagrangian.1} correspond to the non-adhesion conditions. 
Define $\psi : \Re^{m} \to \Re \cup \{ -\infty \}$ by 
\begin{align}
  \psi(\bi{r}) 
  = \inf \{ L(\bi{u},\bi{r}) \mid \bi{u} \in \Re^{d} \} ,
\end{align}
which is the Lagrange dual function. 
The Lagrange dual problem of \eqref{P.contact.1} is then formulated as 
follows: 
\begin{subequations}\label{P.dual.contact.1}%
  \begin{alignat}{3}
    &  \MAX &{\quad}& 
    \psi(\bi{r}) \\
    & \ST && 
    \bi{r} \le \bi{0} . 
  \end{alignat}
\end{subequations}
Since $\psi$ is the pointwise infimum of a family of affine functions of 
$\bi{r}$, it is concave (even if problem \eqref{P.contact.1} is not 
convex). 
Therefore, the dual problem \eqref{P.dual.contact.1} is convex. 

Let $\alpha > 0$ be an arbitrary constant. 
A point $\bi{r} \in \Re^{m}$ is optimal for problem 
\eqref{P.dual.contact.1} if and only if it satisfies 
\begin{align}
  \bi{r} = \Pi_{-}(\bi{r} + \alpha \nabla \psi(\bi{r})) . 
  \label{eq.contact.optimality.1}
\end{align}
This fixed point relation yields the iteration of the projected gradient 
method as follows (see, e.g., \cite{PB14}): 
\begin{align}
  \bi{r}^{(k+1)} 
  := \Pi_{-}(\bi{r}^{(k)} + \alpha \nabla \psi(\bi{r}^{(k)})) . 
  \label{eq.contact.optimality.2}
\end{align}
Here, $\alpha$ plays a role of the step size. 
For given $\bar{\bi{r}} \in \Re_{-}^{m}$, define $\bi{u}_{\bar{\bi{r}}}$ by 
\begin{align}
  \bi{u}_{\bar{\bi{r}}}
  = \argmin \{ L(\bi{u},\bar{\bi{r}}) \mid \bi{u} \in \Re^{d} \} .
\end{align}
It can be shown that the relations 
\begin{align}
  \pdif{}{r_{i}} \psi(\bar{\bi{r}})
  = g_{i}(\bi{u}_{\bar{\bi{r}}}) , 
  \quad i=1,\dots,m
  \label{eq.gradient.contact.1}
\end{align}
holds, when $\pi$ is convex and $g_{1},\dots,g_{m}$ are 
concave \cite[Theorem~9.3-3]{Cia89}. 
Substitution of \eqref{eq.gradient.contact.1} into 
\eqref{eq.contact.optimality.2} results in step~\ref{alg:contact.1.gradient} 
of \refalg{alg:contact.1}. 
Thus, the Uzawa method for problem \eqref{P.contact.1} is viewed as the 
projected gradient method applied to the Lagrange dual problem 
\eqref{P.dual.contact.1}. 
It is worth noting that a reasonable stopping criterion, 
\begin{align*}
  \| \bi{r}^{(k)} - \bi{r}^{(k+1)} \| \le \epsilon
\end{align*}
with threshold $\epsilon$, is derived from the optimality condition in 
\eqref{eq.contact.optimality.1}.


If we assume the small deformation, the total potential energy is given 
by 
\begin{align*}
  \pi(\bi{u}) 
  = \frac{1}{2} \bi{u}^{\top} K \bi{u} - \bi{p}^{\top} \bi{u} . 
\end{align*}
Here, $K \in \Re^{d \times d}$ is the stiffness matrix, which is a 
constant positive definite symmetric matrix, and 
$\bi{p} \in \Re^{d}$ is the external nodal force vector. 
Moreover, $g_{i}$ is linearized as 
\begin{align}
  g_{i}(\bi{u}) = h_{i} - \bi{n}_{i}^{\top} \bi{u} ,
  \quad  i=1,\dots,m. 
  \label{eq.gap.2}
\end{align}
For notational simplicity, we rewrite \eqref{eq.gap.2} as 
\begin{align*}
  \bi{g}(\bi{u})  = \bi{h} - N \bi{u} ,
\end{align*}
where $\bi{h} \in \Re^{m}$ is a constant vector, and 
$N \in \Re^{m \times d}$ is a constant matrix. 
The upshot is that, in the small deformation theory, the frictionless 
contact problem is reduced to the following convex quadratic programming 
(QP) problem: 
\begin{subequations}\label{P.contact.2}%
  \begin{alignat}{3}
    &  \MIN &{\quad}& 
    \frac{1}{2} \bi{u}^{\top} K \bi{u} - \bi{p}^{\top} \bi{u} \\
    & \ST && 
    \bi{h} - N \bi{u}  \ge \bi{0} . 
  \end{alignat}
\end{subequations}
The step size, $\alpha$, of the Uzawa method for solving problem 
\eqref{P.contact.2} is chosen as follows. 
Let $\lambda_{1}(K)$ and $\sigma_{d}(N)$ denote the minimum eigenvalue 
of $K$ and the maximum singular value of $\sigma_{d}(N)$. 
If $\alpha$ is chosen so that 
\begin{align*}
  \alpha \in ] 0, 2 \lambda_{1}(K) / \sigma_{d}(N) [ , 
\end{align*}
then it is shown that the Uzawa method converges \cite[section~9.4]{Cia89}.

\section{Accelerated Uzawa method}
\label{sec:acceleration}

In this section, we present an accelerated version of the Uzawa method. 

The Uzawa method in \refalg{alg:contact.1} is essentially viewed as the 
projected gradient method \eqref{eq.contact.optimality.2}. 
This method is considered a special case of the proximal gradient method, 
because the proximal operator of the indicator function of a nonempty 
closed convex set is reduced to the projection onto the set \citep{PB14}. 
Therefore, it is natural to apply the acceleration scheme for the 
proximal gradient method \citep{BT09} to the projected gradient update 
\eqref{eq.contact.optimality.2}. 
This results in the update 
\begin{align}
  \bi{r}^{(k+1)} 
  &:= \Pi_{-}(\bi{\rho}^{(k)} + \alpha \nabla \psi(\bi{\rho}^{(k)})) , 
  \label{eq.accelerated.gradient.1.1} \\
  \bi{\rho}^{(k+1)} 
  &:= \bi{r}^{(k)} 
  + \omega_{k+1}(\bi{r}^{(k+1)}  - \bi{r}^{(k)}) . 
  \label{eq.accelerated.gradient.1.2}
\end{align}
Here, $\omega_{k} \in [0,1)$ is an extrapolation parameter which is to 
be determined so that the convergence acceleration is achieved. 
In \cite{BT09}, $\omega_{k}$ is chosen as 
\begin{align*}
  \tau_{k+1} 
  &:= \frac{1}{2} \Bigl(
  1 + \sqrt{1 + \tau_{k}^{2}}  \Bigr) , \\
  \omega_{k+1} 
  &:= \frac{\tau_{k} - 1}{\tau_{k+1}}  
\end{align*}
with $\tau_{0}:=1$. 
The sequence of the dual objective values, $\{ \psi(\bi{r}^{(k)}) \}$, 
converges to the optimal value with rate $O(1/k^{2})$ \citep{BT09}. 



The accelerated method in \eqref{eq.accelerated.gradient.1.1} and 
\eqref{eq.accelerated.gradient.1.2} is not guaranteed to be monotone in 
the dual objective value. 
Therefore, we follow \citet{OdC15} in incorporating an adaptive restart 
technique of the acceleration scheme. 
Namely, we perform the restart procedure whenever the momentum term, 
$\bi{r}^{(k+1)} - \bi{r}^{(k)}$, and the gradient of the objective 
function, $\nabla\psi(\bi{r}^{(k)})$, make an obtuse angle, i.e., 
\begin{align*}
  \nabla\psi(\bi{r}^{(k)})^{\top} 
  (\bi{r}^{(k+1)} - \bi{r}^{(k)}) < 0 . 
\end{align*}
As the upshot, the accelerated Uzawa method with adaptive restart is 
listed in \refalg{alg:contact.3}. 
One reasonable stopping criterion is 
\begin{align}
  \| \bi{\rho}^{(k)} - \bi{r}^{(k+1)} \| \le \epsilon  
  \label{eq.stopping.accelerate}
\end{align}
at step~\ref{alg:contact.3.reaction}.

\begin{algorithm}
  \caption{accelerated Uzawa method with restart}
  \label{alg:contact.3}
  \begin{algorithmic}[1]
    \Require
    $\alpha > 0$, $\bi{r}^{(0)} \in \Re_{-}^{m}$, 
    $\bi{\rho}^{(0)} := \bi{r}^{(0)}$, 
    $\tau_{0} := 1$
    \For{$k=0,1,\dots$}
    \State \label{alg:contact.3.equilibrium}
    $\bi{u}^{(k)}$ solves 
    $\displaystyle 
    \nabla \pi(\bi{u}) + \sum_{i=1}^{m} \rho_{i}^{(k)}\nabla g_{i}(\bi{u})  = \bi{0}$
    \State \label{alg:contact.3.gap}
    $\bi{\gamma}^{(k)} := \bi{g}(\bi{u}^{(k)})$
    \State \label{alg:contact.3.reaction}
    $\bi{r}^{(k+1)} 
    := \min \{ \bi{0}, \bi{\rho}^{(k)} + \alpha \bi{\gamma}^{(k)} \}$ 
    \State 
    $\displaystyle
    \tau_{k+1} := \frac{1}{2}
    \Bigl(1 + \sqrt{1 + 4\tau_{k}^{2}} \Bigr)$
    \If{$(\bi{\gamma}^{(k)})^{\top} (\bi{r}^{(k+1)} - \bi{r}^{(k)}) \ge 0$} 
    \State 
    $\displaystyle
    \bi{\rho}^{(k+1)} 
    := \bi{r}^{(k+1)} + \frac{\tau_{k}-1}{\tau_{k+1}} 
    (\bi{r}^{(k+1)} - \bi{r}^{(k)})$ 
    \Else
    \State
    $\bi{\rho}^{(k+1)} := \bi{r}^{(k+1)}$ 
    \State
    $\tau_{k+1} := 1$
    \EndIf
    \EndFor
  \end{algorithmic}
\end{algorithm}

In the case of the small deformation theory, 
steps~\ref{alg:contact.3.equilibrium} and \ref{alg:contact.3.gap} 
of \refalg{alg:contact.3} are simplified as follows. 
\begin{itemize}
  \item Step~\ref{alg:contact.3.equilibrium}: 
        Let $\bi{u}^{(k)}$ be the solution to the system of linear 
        equations 
\begin{align}
  K \bi{u} = \bi{p} + N^{\top} \bi{\rho}^{(k)} . 
  \label{eq.linear.equation}
\end{align}
        Here, the coefficient matrix $K$ is common to all the iterations. 
        Hence, we carry out the Cholesky factorization only at the first 
        iteration; at the following iterations, we can solve 
        \eqref{eq.linear.equation} only with the back-substitutions. 
  \item Step~\ref{alg:contact.3.gap}: 
        $\bi{\gamma}^{(k)} := \bi{h} - N \bi{u}^{(k)}$. 
\end{itemize}

\section{Preliminary numerical experiments}
\label{sec:ex}

\begin{figure}[tp]
  \centering
  \includegraphics[scale=0.40]{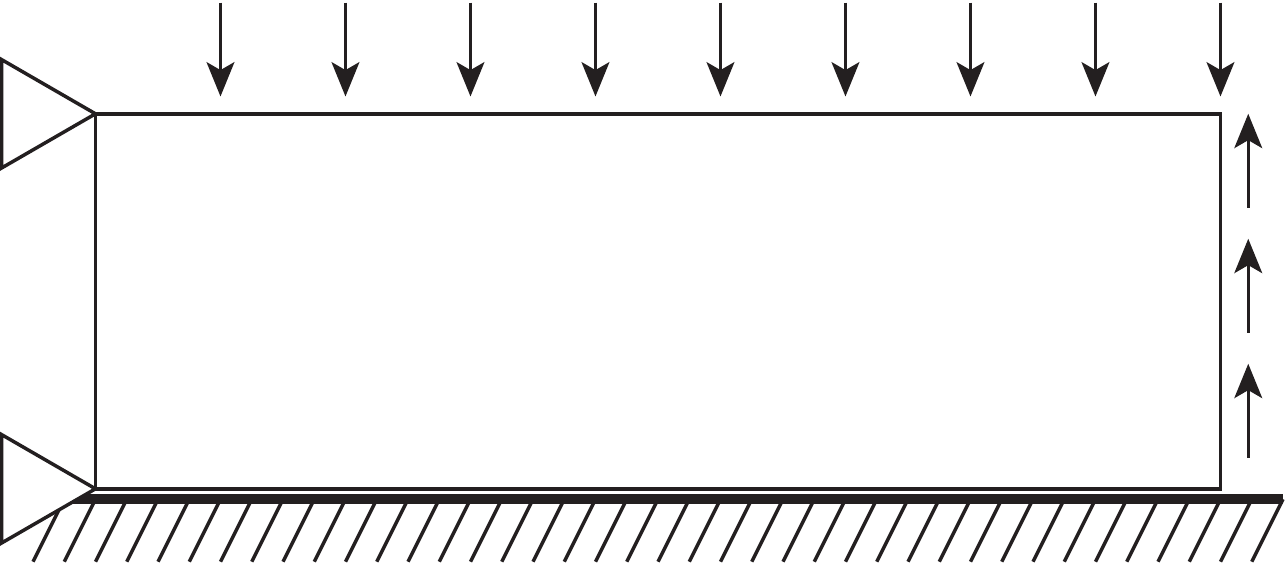}
  \caption{An elastic body on the obstacle.}
  \label{fig:ex_body}
\end{figure}

\begin{figure}[tp]
  \centering
  \subfigure[]{\label{fig:convergence_dual_obj}
  \includegraphics[scale=0.50]{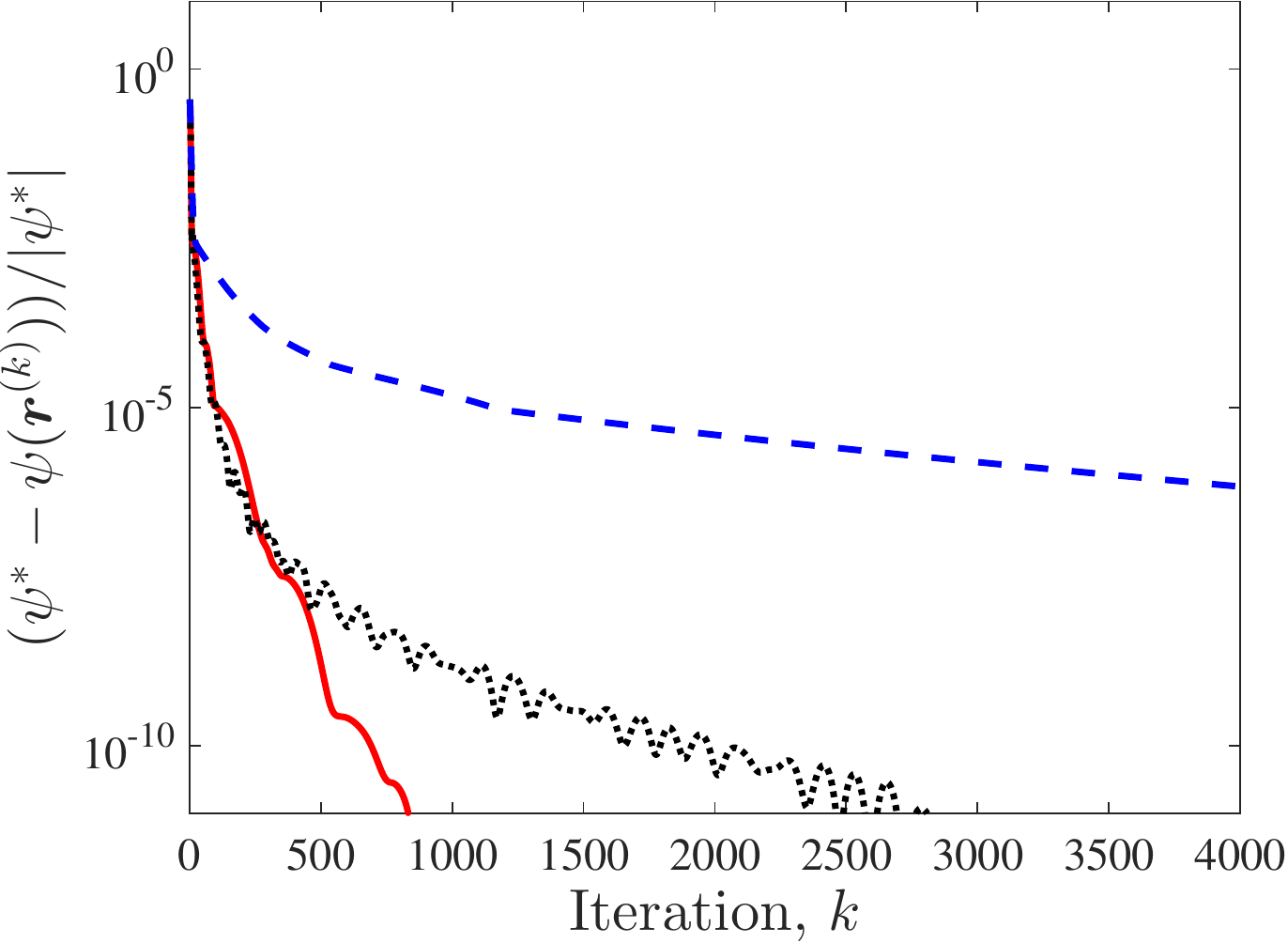}
  }
  \hspace{1em}
  \subfigure[]{\label{fig:convergence_residual}
  \includegraphics[scale=0.50]{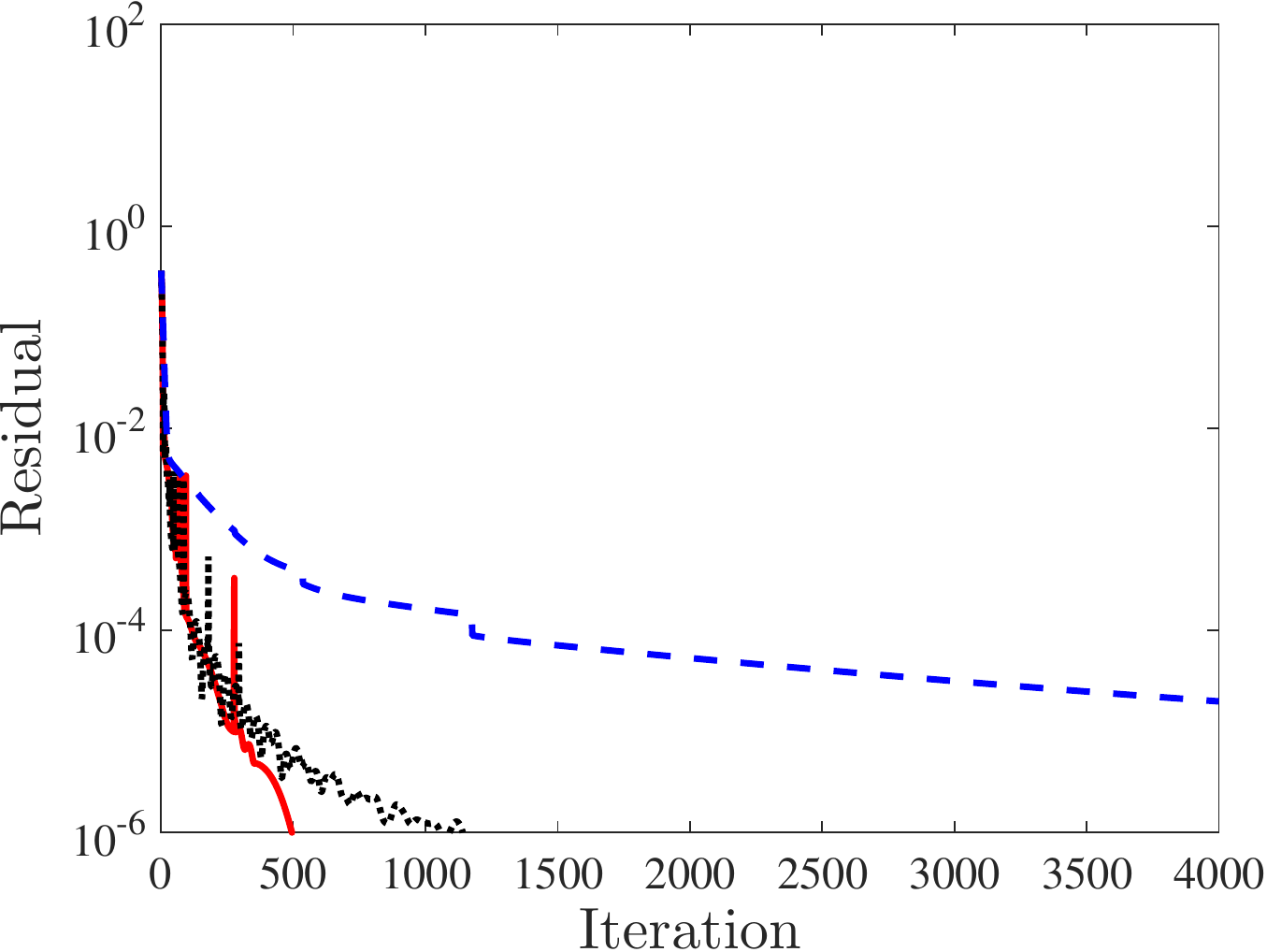}
  }
  \caption{Convergence history for $(N_{X},N_{Y})=(30,10)$. 
  \subref{fig:convergence_dual_obj} The dual objective value; and 
  \subref{fig:convergence_residual} the total residual. 
  ``{\em Solid line\/}'' the accelerated Uzawa method with restart; 
  ``{\em dotted line\/}'' the accelerated Uzawa method without restart; and 
  ``{\em dashed line\/}'' the Uzawa method. }
  \label{fig:convergence}
\end{figure}

Consider an elastic body shown in \reffig{fig:ex_body}. 
The body is in the plane-stress state, with 
thickness $5\,\mathrm{mm}$, width $60\,\mathrm{mm}$, and 
height $20\,\mathrm{mm}$. 
It consists of an isotropic homogeneous material with Young's modulus 
$200\,\mathrm{GPa}$ and Poisson's ratio $0.3$. 
The bottom edge of the body is on the rigid obstacle, and the left 
edge is fixed by supports. 
The uniform downward traction of $50\,\mathrm{kPa}$ is applied to the 
top edge, and the uniform upward traction of $500\,\mathrm{kPa}$ is 
applied to the right edge. 
The body is discretized into $N_{X} \times N_{Y}$ four-node 
quadrilateral (Q4) finite elements, where $N_{X}$ $(=3N_{Y})$ is varied to 
generate problem instances with different sizes. 
The number of degrees of freedom of displacements is 
$d=2N_{X}(N_{Y}+1)$ and 
the number of contact candidate nodes is $m=N_{X}$. 
At the equilibrium state, about 73.3\% of contact candidate nodes are 
in contact with nonzero reactions. 
We assume the small deformation, and solve QP problem \eqref{P.contact.2}. 
Computation was carried out on a $2.2\,\mathrm{GHz}$ Intel Core i5 
processor with $8\,\mathrm{GB}$ RAM. 

The proposed algorithm was implemented with MATLAB ver.\ 9.0.0. 
The threshold for the stopping criteria in 
\eqref{eq.stopping.accelerate} is set to $\epsilon=10^{-6}$. 
Comparison was performed with \textsf{QUADPROG} \citep{MATLAB}, 
\textsf{SDPT3} ver.\ 4.0 \citep{TTT03} 
with \textsf{CVX} ver.\ 2.1 \citep{GB08}, and 
\textsf{PATH} ver.\ 4.7.03 \citep{FM00}
via the MATLAB Interface \citep{FM99}. 
\textsf{QUADPROG} is a MATLAB built-in function for QP. 
We use an implementation of an interior-point method 
with setting the termination threshold, the parameter \texttt{TolFun}, 
to $10^{-8}$. 
\textsf{SDPT3} implements a primal-dual interior-point method for
solving conic programming problems. 
A QP problem is converted to the standard form of conic programming 
by \textsf{CVX}, where the parameter \texttt{cvx\_precision} for 
controlling the solver precision is set to \texttt{high}. 
\textsf{PATH} is a nonsmooth Newton method to solve mixed 
complementarity problems, and hence can solve the KKT condition for QP. 

\begin{figure}[tp]
  \centering
  \subfigure[]{\label{fig:time_all_compare}
  \includegraphics[scale=0.50]{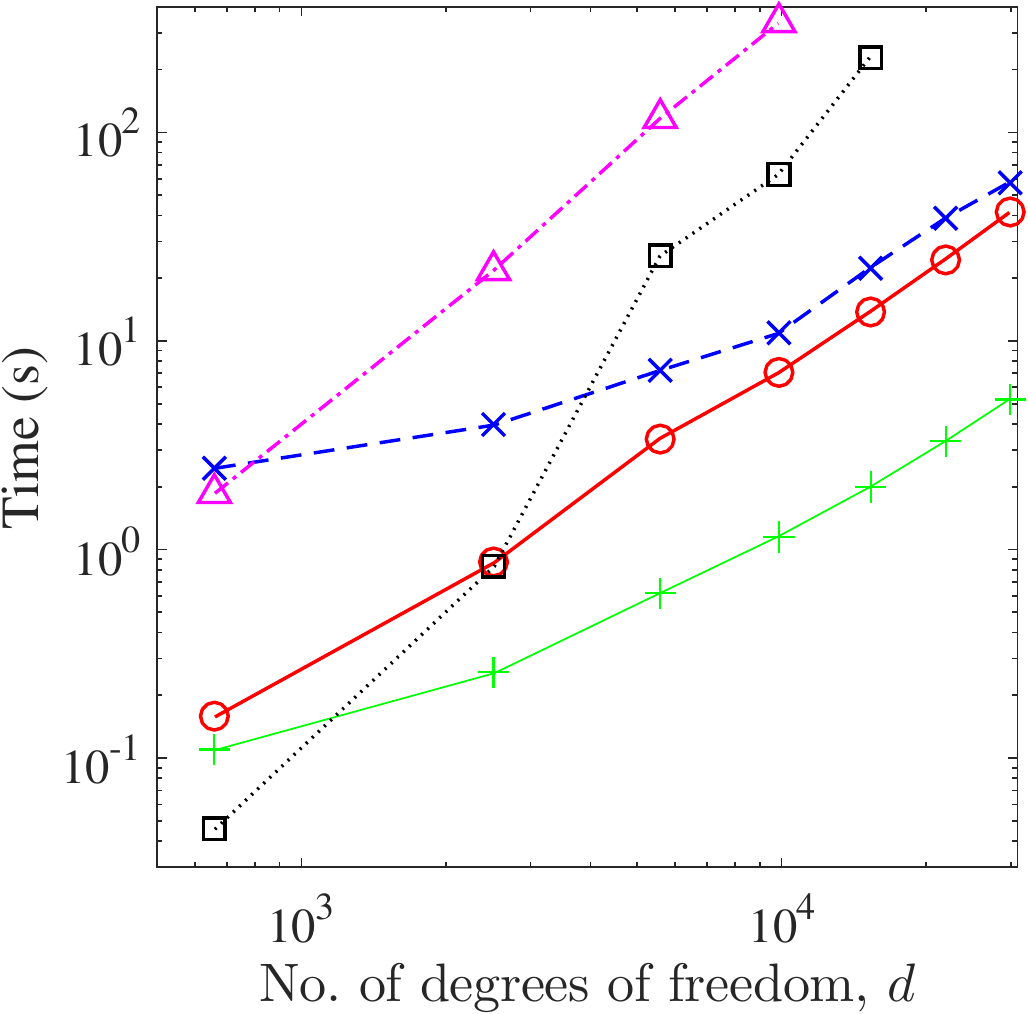}
  }
  \hspace{1em}
  \subfigure[]{\label{fig:residual_all_compare}
  \includegraphics[scale=0.50]{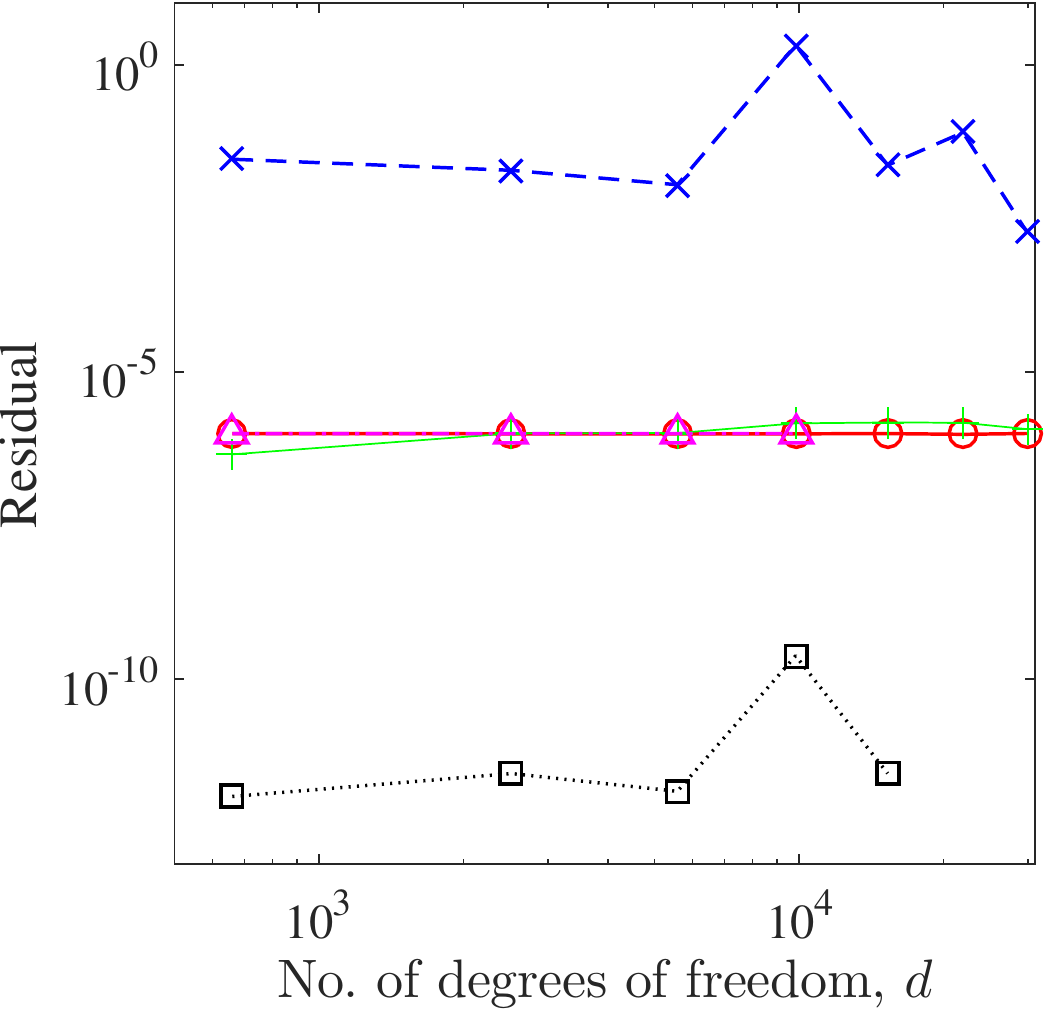}
  }
  \caption{Comparison of 
  \subref{fig:time_all_compare} the computational time; and 
  \subref{fig:residual_all_compare} the total residual. 
  ``$\circ$'' The accelerated Uzawa method with restart; 
  ``$\triangle$'' the Uzawa method; 
  ``$\times$'' \textsf{SDPT3}; 
  ``$+$'' \textsf{QUADPROG}; and 
  ``$\square$'' \textsf{PATH}. 
  }
  \label{fig:comparison}
\end{figure}

\reffig{fig:convergence} reports the convergence history of the proposed 
algorithm (\refalg{alg:contact.3}). 
It also shows the result of the accelerated Uzawa method without restart 
scheme, and that of \refalg{alg:contact.1} (i.e., the Uzawa method 
without acceleration). 
\reffig{fig:convergence_dual_obj} shows the convergence history of the 
dual objective function. 
It is observed that the acceleration and restart schemes drastically speed
up the convergence. 
Particularly, with the proposed algorithm, the dual objective value 
seems to converge quadratically and monotonically. 
The residual in \reffig{fig:convergence_residual} is defined by using 
the KKT condition. 
Namely, we define vectors $\bi{e}_{j}$ $(j=1,2,3,4)$ by 
\begin{align*}
  \bi{e}_{1} 
  &= K \bi{u}^{(k)} - \bi{f} - N^{\top} \bi{r}^{(k)} , \\
  \bi{e}_{2} 
  &= \min \{ \bi{g}(\bi{u}^{(k)}), \bi{0} \} , \\
  \bi{e}_{3} 
  &= \max \{ \bi{r}^{(k)}, \bi{0} \} , \\
  \bi{e}_{4} 
  &= -\diag(\bi{g}(\bi{u}^{k})) \bi{r}^{(k)} , 
\end{align*}
and the value of $\| (\bi{e}_{1},\bi{e}_{2},\bi{e}_{3},\bi{e}_{4}) \|$ 
is reported in \reffig{fig:convergence_residual}. 

\reffig{fig:comparison}  compares the computational results of the five 
methods. 
It is observed in \reffig{fig:time_all_compare} that \textsf{QUADPROG} 
is fastest for almost all the problem instances. 
The computational time required by \textsf{PATH} is very small for small 
instances, but drastically increases as the instance size increases. 
On the other hand, it is observed in \reffig{fig:residual_all_compare} 
that \textsf{PATH} converges to the highest accuracy. 
For large instances, the proposed method is the second-fastest method 
among the five methods. 
The computational time of \textsf{SDPT3} is larger than that of the 
proposed method. 
Also, \textsf{SDPT3} converges to the lowest accuracy. 
The residuals of the solutions obtained by the proposed method, the 
unaccelerated Uzawa method, and \textsf{QUADPROG} are similar. 
The computational time required by the unaccelerated Uzawa method is 
very large. 

The upshot is that, within the small deformation theory, the proposed 
method is quite efficient, but the interior-point method for QP 
(\textsf{QUADPROG}) is more efficient. 
If we consider large deformation, problem \eqref{P.contact.1} is not 
convex, and hence \textsf{QUADPROG} and \textsf{SDPT3} cannot be adopted. 
Efficiency of the proposed method, that solves the convex dual problem, 
remains to be studied.

\section{Concluding remarks}
\label{sec:conclude}

This paper has presented an accelerated Uzawa method. 
The method can be recognized as an accelerated projected gradient method 
solving the Lagrangian dual problem of a constrained optimization 
problem. 
It has been shown in the numerical experiments that the acceleration and 
restart schemes drastically speed up the convergence of the Uzawa method. 
Besides, the proposed method is easy to implement. 

Many possibilities of extensions could be considered. 
To update the primal variables, the Uzawa method solves a system of 
linear equations. 
For large-scale problems, this process might be replaced by a fast 
first-order minimization method of a convex quadratic function, 
e.g., \cite{LS13}. 
Also, an extension to large deformation problems remains to be explored.

\bigskip
\paragraph{Acknowledgments}

This work is partially supported by 
JSPS KAKENHI 26420545 and 17K06633.


\begin{thebibliography}{99}

\bibitem[\protect\citeauthoryear{Alart and Curnier}{1991}]{AC91}
  {P.~Alart, A.~Curnier}:
  {A mixed formulation for frictional contact problems prone to Newton
    like solution methods}.
  {\em Computer Methods in Applied Mechanics and Engineering},
  \textbf{92}, 353--375 (1991).

\bibitem[\protect\citeauthoryear{Allaire}{2007}]{All07}
  {G.~Allaire}:
  {\em Numerical Analysis and Optimization}.
  Oxford University Press, New York (2007).

\bibitem[\protect\citeauthoryear{Beck and Teboulle}{2009}]{BT09}
  {A.~Beck, M.~Teboulle}:
  {A fast iterative shrinkage-thresholding algorithm 
    for linear inverse problems}.
  {\em SIAM Journal on Imaging Sciences},
  \textbf{2}, 183--202 (2009).

\bibitem[\protect\citeauthoryear{Becker {\em et al.\/}}{2011}]{BBC11}
  {S.~Becker, J.~Bobin, E.~J.~Cand\`{e}s}:
  {NESTA: A fast and accurate first-order method for sparse recovery}.
  {\em SIAM Journal on Imaging Sciences},
  \textbf{4}, 1--39 (2011).



\bibitem[\protect\citeauthoryear{Ciarlet}{1989}]{Cia89}
  {P.~G.~Ciarlet}:
  {\em Introduction to Numerical Linear Algebra and Optimisation}.
  Cambridge University Press, Cambridge (1989).



\bibitem[\protect\citeauthoryear{Duvaut and Lions}{1976}]{DL76}
  {G.~Duvaut, J.~L.~Lions}:
  {\em Inequalities in Mechanics and Physics}.
  Springer-Verlag, Berlin (1976).

\bibitem[\protect\citeauthoryear{Ferris and Munson}{1999}]{FM99}
  {M.~C.~Ferris, T.~S.~Munson}:
  {Interfaces to PATH 3.0: Design, implementation and usage}.
  {\em Computational Optimization and Applications},
  \textbf{12}, 207--227 (1999).

\bibitem[\protect\citeauthoryear{Ferris and Munson}{2000}]{FM00}
  {M.~C.~Ferris, T.~S.~Munson}:
  {Complementarity problems in GAMS and the PATH solver}.
  {\em Journal of Economic Dynamics and Control},
  \textbf{24}, 165--188 (2000).

\bibitem[\protect\citeauthoryear{Grant and Boyd}{2008}]{GB08}
  {M.~Grant, S.~Boyd}:
  {Graph implementations for nonsmooth convex programs}.
  In: V.~Blondel, S.~Boyd, H.~Kimura (eds.),
  {\em Recent Advances in Learning and Control (A Tribute to M.~Vidyasagar)},
  Springer (2008). pp.~95--110.


\bibitem[\protect\citeauthoryear{Haslinger {\em et al.\/}}{2004}]{HKD04}
  {J.~Haslinger, R.~Ku\v{c}era, Z.~Dost\'{a}l}:
  {An algorithm for the numerical realization of 3D contact 
    problems with Coulomb friction}.
  {\em Journal of Computational and Applied Mathematics},
   \textbf{164--165}, 387--408 (2004).

\bibitem[\protect\citeauthoryear{Kanno}{2016}]{Kan16}
  {Y.~Kanno}:
  {A fast first-order optimization approach to elastoplastic analysis 
      of skeletal structures}.
  {\em Optimization and Engineering},
  \textbf{17}, 861--896 (2016).


\bibitem[\protect\citeauthoryear{Koko}{2009}]{Kok09}
  {J.~Koko}:
  {Uzawa block relaxation domain decomposition method for a two-body
    frictionless contact problem}.
  {\em Applied Mathematics Letters},
  \textbf{22}, 1534--1538 (2009).



\bibitem[\protect\citeauthoryear{Lee and Sidford}{2013}]{LS13}
  {Y.~T.~Lee, A.~Sidford}:
  {Efficient accelerated coordinate descent methods and 
    faster algorithms for solving linear systems}.
  {\em IEEE 54th Annual Symposium on Foundations of Computer Science}, 
    Berkeley (2013) pp.~147--156.

\bibitem[\protect\citeauthoryear{Nesterov}{1983}]{Nes83}
  {Y.~Nesterov}:
  {A method of solving a convex programming problem with 
    convergence rate $O(1/k^{2})$}.
  {\em Soviet Mathematics Doklady},
  \textbf{27}, 372--376 (1983).

\bibitem[\protect\citeauthoryear{Nesterov}{2004}]{Nes04}
  {Y.~Nesterov}:
  {\em Introductory Lectures on Convex Optimization: A Basic Course}.
  Kluwer Academic Publishers, Dordrecht (2004).

\bibitem[\protect\citeauthoryear{O'Donoghue and Cand\`{e}s}{2015}]{OdC15}
  {B.~O'Donoghue, E.~Cand\`{e}s}:
  {Adaptive restart for accelerated gradient schemes}.
  {\em Foundations of Computational Mathematics},
  \textbf{15}, 715--732 (2015).

\bibitem[\protect\citeauthoryear{Parikh and Boyd}{2014}]{PB14}
  {N.~Parikh, S.~Boyd}:
  {Proximal algorithms}.
  {\em Foundations and Trends in Optimization},
  \textbf{1}, 127--239 (2014).

\bibitem[\protect\citeauthoryear{Puso {\em et al.\/}}{2008}]{PLS08}
  {M.~A.~Puso, T.~A.~Laursen, J.~Solberg}:
  {A segment-to-segment mortar contact method for quadratic 
    elements and large deformations}.
  {\em Computer Methods in Applied Mechanics and Engineering},
  \textbf{197}, 555--566 (2008).

\bibitem[\protect\citeauthoryear{Raous}{1999}]{Rao99}
  {M.~Raous}:
  {Quasistatic Signorini problem with Coulomb friction and 
    coupling to adhesion}.
  In: P.~Wriggers, P.~Panagiotopoulos (eds.),
  {\em New Developments in Contact Problems},
  Springer-Verlag, Wien (1999) pp.~101--178.

\bibitem[\protect\citeauthoryear{Rudoy}{2015}]{Rud15}
  {E.~M.~Rudoy}:
  {Domain decomposition method for a model crack problem
    with a possible contact of crack edges}.
  {\em Computational Mathematics and Mathematical Physics},
  \textbf{55}, 305--316 (2015).

\bibitem[\protect\citeauthoryear{Temizer}{2012}]{Tem12}
  {\.{I}.~Temizer}:
  {A mixed formulation of mortar-based frictionless contact}.
  {\em Computer Methods in Applied Mechanics and Engineering},
  \textbf{223--224}, 173--185 (2012).

\bibitem[\protect\citeauthoryear{Temizer {\em et al.\/}}{2012}]{TWH12}
  {\.{I}.~Temizer, P.~Wriggers, T.~J.~R.~Hughes}:
  {Three-dimensional mortar-based frictional contact treatment 
    in isogeometric analysis with NURBS}.
  {\em Computer Methods in Applied Mechanics and Engineering},
  \textbf{209--212}, 115--128 (2012).

\bibitem[\protect\citeauthoryear{The MathWorks, Inc.}{2016}]{MATLAB}
  {The MathWorks, Inc.}:
  {\em MATLAB Documentation}.
  \url{http://www.mathworks.com/}
  (Accessed November 2016).

\bibitem[\protect\citeauthoryear{T\"{u}t\"{u}nc\"{u} {\em et al.\/}}{2003}]{TTT03}
  {R.~H.~T\"{u}t\"{u}nc\"{u}, K.~C.~Toh, M.~J.~Todd}:
  {Solving semidefinite-quadratic-linear programs using SDPT3}.
  {\em Mathematical Programming},
  \textbf{B95}, 189--217 (2003).

\bibitem[\protect\citeauthoryear{Uzawa}{1958}]{Uza58}
  {H.~Uzawa}:
  {Iterative methods for concave programming}.
  In: K.~J.~Arrow, L.~Hurwicz, H.~Uzawa (eds.),
  {\em Studies in Linear and Non-Linear Programming},
  Stanford University Press, Stanford (1958) pp.~154--165.

\bibitem[\protect\citeauthoryear{Vikhtenko and Namm}{2007}]{VN07}
  {E.~M.~Vikhtenko, R.~V.~Namm}:
  {Duality scheme for solving the semicoercive Signorini
    problem with friction}.
  {\em Computational Mathematics and Mathematical Physics},
  \textbf{47}, 1938--1951 (2007).

\bibitem[\protect\citeauthoryear{Wriggers}{2006}]{Wri06}
  {P.~Wriggers}:
  {\em Computational Contact Mechacnics (2nd ed.)}.
  Springer-Verlag, Berlin (2006).

\bibitem[\protect\citeauthoryear{Zavarise and De Lorenzis}{2012}]{ZDl12}
  {G.~Zavarise, L.~De Lorenzis}:
  {An augmented Lagrangian algorithm for contact mechanics based 
    on linear regression}.
  {\em International Journal for Numerical Methods in Engineering},
  \textbf{91}, 825--842 (2012).

\end{thebibliography}
\end{document}